\documentclass[11pt]{article}
\usepackage{amsmath,amssymb,amscd,array}

\let\ssection=\section
\renewcommand{\section}{\setcounter{equation}{0}\ssection}

\def\d{\delta}

\def\om{\omega}
\def\r{\rho}
\def\a{\alpha}
\def\b{\beta}
\def\s{\sigma}
\def\vfi{\varphi}

\def\l{\lambda}
\def\m{\mu}

\def\implies{\Rightarrow}

\newcommand{\bbR}{\mathbb{R}}

\newcommand{\Diff}{\mathrm{Diff}}

\newcommand{\cF}{{\mathcal{F}}}
\newcommand{\cD}{{\mathcal{D}}}

\newcommand{\Hom}{\mathrm{Hom}}

\newcommand{\SL}{\mathrm{SL}}
\newcommand{\Sl}{\mathrm{sl}}

\newcommand{\Vect}{\mathrm{Vect}}

\begin{document}

\frenchspacing

\def\d{\delta}
\def\g{\gamma}
\def\om{\omega}
\def\r{\rho}
\def\a{\alpha}
\def\b{\beta}
\def\s{\sigma}
\def\vfi{\varphi}
\def\l{\lambda}
\def\m{\mu}
\def\implies{\Rightarrow}

\oddsidemargin .1truein
\newtheorem{thm}{Theorem}[section]
\newtheorem{lem}[thm]{Lemma}
\newtheorem{cor}[thm]{Corollary}
\newtheorem{pro}[thm]{Proposition}
\newtheorem{ex}[thm]{Example}
\newtheorem{rmk}[thm]{Remark}
\newtheorem{defi}[thm]{Definition}
\title{Cohomology of the vector fields Lie algebras on $\mathbb{RP}^1$
acting on bilinear differential operators}
\author{Sofiane Bouarroudj\\
{\footnotesize Department of Mathematics, U.A.E. University,
Faculty of Science} \\
{\footnotesize P.O. Box 15551, Al-Ain, United Arab Emirates.}\\
{\footnotesize  e-mail:bouarroudj.sofiane@uaeu.ac.ae} }
\date{}
\maketitle
\begin{abstract} The main topic of this paper is two folds.
First, we compute the first relative cohomology group of the Lie
algebra of smooth vector fields on the projective line,
$\Vect(\mathbb{RP}^1),$ with coefficients in the space of bilinear
differential operators that act on tensor densities, ${\cal
D}_{\lambda, \nu;\mu},$ vanishing on the Lie algebra
$\Sl(2,\mathbb{R}).$ Second, we compute the first cohomology group
of the Lie algebra $\Sl(2,\mathbb{R})$ with coefficients in ${\cal
D}_{\lambda, \nu;\mu}.$
\end{abstract}
\section{Introduction}
Let ${\mathfrak g}$ be a Lie algebra and let ${\cal M}$ and ${\cal
N}$ be two ${\mathfrak g}$-modules. It is well-known that
nontrivial extensions of ${\mathfrak g}$-modules:
$$
0\rightarrow {\cal
M}\rightarrow \cdot \rightarrow {\cal N}\rightarrow 0
$$
are classified by the first cohomology group ${\mathrm
H}^1({\mathfrak g}; \Hom ({\cal N},{\cal M}))$ (see, e.g.,
\cite{f}). Any 1-cocycle ${\cal L}$ generates a new action on
${\cal M}\oplus {\cal N}$ as follows: for all $g\in {\mathfrak g}$
and for all $(a,b)\in {\cal M}\oplus {\cal N},$ we define $
g^*(a,b):=(g^*a+ C^{st} {\cal L}(b), g^*b).$ For the space of
tensor densities of weight $\lambda,$ ${\cal F}_\lambda,$ viewed
as a module over the Lie algebra of smooth vector fields
$\Vect(\mathbb{RP}^1),$ the classification of nontrivial
extensions
$$
0\rightarrow {\cal F}_{\mu}\rightarrow \cdot \rightarrow {\cal
F}_{\lambda}\rightarrow 0,
$$
leads Feigin and Fuks in \cite{gf} to compute the cohomology group
${\mathrm H}^1(\Vect(\mathbb{RP}^1); \Hom ({\cal
F}_{\lambda},{\cal F}_{\mu})).$ Later, Ovsienko and the author in
\cite{bo1} have computed the corresponding relative cohomology
group with respect to $\Sl(2,\mathbb{R}),$ namely
\begin{equation}
\label{chmes} {\mathrm H}^1(\Vect(\mathbb{RP}^1),
\Sl(2,\mathbb{R}); \Hom_{\mathrm {diff} } ({\cal
F}_{\lambda},{\cal F}_{\mu})).
\end{equation}
In fact, the study of the cohomology group (\ref{chmes}) has
arisen from the study of the equivariant quantization procedure
introduced in \cite{cmz, lol}. It has been proved that there
exists an $\Sl(2,\mathbb{R})$-equivariant quantization map from
the space of symbols and to the space of differential operators,
but it is not $\Vect(\mathbb{RP}^1)$-equivariant. The obstruction
here is given by the 1-cocycles that span the cohomology group
(\ref{chmes}) (see \cite{bo1, gar}). The computation is based on
an old result of Gordan \cite{g} on the classification of
$\Sl(2,\mathbb{R})$-invariant bilinear differential operators that
act on tensor densities. Moreover, the case of a
higher-dimensional manifold has been studied in \cite{b, lo2}, and
the case of a Riemann surface has been studied in \cite{bg}. In
this paper, we will compute the first cohomology group
$$
{\mathrm H}^1(\Vect(\mathbb{RP}^1), \Sl(2,\mathbb{R});
\Hom_{\mathrm{diff}} ({\cal F}_{\lambda}\otimes {\cal
F}_{\nu},{\cal F}_{\mu})).
$$
It turns out that the dimension of the cohomology group above can
reach three, for some particular values of $\lambda$ and $\nu,$
which is a contradistinction with that of the cohomology group
(\ref{chmes}) in which the dimension is almost one.

Moreover, we compute the cohomology group
$$
{\mathrm H}^1(\Sl(2,\mathbb{R}); \Hom_{\mathrm{diff}} ({\cal
F}_{\lambda}\otimes {\cal F}_{\nu},{\cal F}_{\mu})).
$$
For linear differential operators, the analogue of the cohomology
group above has been studied by Lecomte in \cite{l}.

\section{$\Vect(\mathbb{R})$-module structures on the space of
bilinear differential operators}
Consider the standard action of $\SL(2,\mathbb{R})$ on
$\mathbb{RP}^1$ by projective transformations. It is given in
homogenous coordinates by
$$
x\mapsto \frac{ax+b}{cx+d}, \quad \mbox{ where } \left(
\begin{array}{cc}
a & b\\
c & d
\end{array}
\right )\in \SL(2,\mathbb{R}).
$$
This action generates global vector fields
$$
\frac{d}{dx},\quad x\frac{d}{dx},\quad x^2\frac{d}{dx},
$$ that form a Lie subalgebra of $\Vect(\mathbb{RP}^1),$ isomorphic to
$\Sl(2,\mathbb{R})$ (see e.g. \cite{olv}). Throughout this paper,
$\Sl(2,\mathbb{R})$ will be refered to this subalgebra.
\subsection{The space of tensor densities on $\mathbb{RP}^1$}
The space of tensor densities of weight $\lambda$ on
$\mathbb{RP}^1$, denoted by ${\cal F}_{\l}$, is the space of
sections of the line bundle $ (T^*\mathbb{RP}^1)^{\otimes \l}.$
This space coincides with the space of functions and differential
forms for $\lambda=0$ and for $\lambda=1,$ respectively. The Lie
algebra $\Vect(\mathbb{RP}^1)$ acts on $\cal F_\l$ by the Lie
derivative. For all $X\in \Vect(\mathbb{RP}^1)$ and for all $\phi
\in \cal F_\l:$
\begin{equation}
\label{dens} L_{X}^{\l}(\phi)=X\phi' +\l \phi \,X',
\end{equation}
where the superscript $'$ stands for $d/dx.$
\subsection{The space of bilinear differential
operators as a $\Vect(\mathbb{RP}^1)$-module}
We are interested in defining a three-parameter family of
$\Vect(\mathbb{RP}^1)$-modules on the space of bilinear
differential operators. The counterpart
$\Vect(\mathbb{RP}^1)$-modules of the space of linear differential
operators is a classical object (see e.g. \cite{w}).

Consider bilinear differential operators that act on tensor
densities:
\begin{equation}
\label{Op} A:{\cal F}_{\lambda}\otimes {\cal F}_{\nu}\to {\cal
F}_{\mu}.
\end{equation}
The Lie algebra $\Vect(\mathbb{RP}^1)$ acts on the space of
bilinear differential operators as follows. For all $\phi \in
\cF_\lambda$ and for all $\psi \in \cF_\nu:$
\begin{equation}
\label{act} L_X^{\l,\nu;\mu}(A)(\phi, \psi)=L_X^\mu \circ A(\phi,
\psi) -A (L_X^\l\,\phi,\psi)-A (\phi,L_X^\nu\,\psi).
\end{equation}
where $L_X^\l$ is the action (\ref{dens}). We denote by ${\cal
D}_{\lambda,\mu;\nu}$ the space of bilinear differential operators
(\ref{Op}) endowed with the defined $\Vect(\mathbb{RP}^1)$-module
structure (\ref{act}).
\section{Cohomology of $\Sl(2,\bbR)$ acting on
${\cal D}_{\lambda,\nu;\mu}$}
In this section, we will compute the ``differentiable'' cohomology
of the Lie algebra $\Sl(2,\bbR)$ with coefficients in the space of
bilinear differential operators ${\cal D}_{\lambda,\nu;\mu}.$
Namely, we consider only cochains that are given by differentiable
maps.
\begin{thm}
\label{leco} (i) If $\mu=\lambda+\nu$ then
\begin{equation} \nonumber \mathrm
H^1(\Sl(2,\bbR);  \cD_{\lambda,\nu;\mu})=\mathbb{R}.
\end{equation}
(ii) If $\mu-\lambda-\nu=k,$  where $k$ is a positive integer,
then
\begin{equation} \nonumber \mathrm
H^1(\Sl(2,\bbR);  \cD_{\lambda,\nu;\mu})= \left\{
\begin{array}{ll}
\bbR^3,& \mbox{if }\, (\lambda,\nu)=(-\frac{s}{2},-\frac{t}{2}),
\mbox{where }  0\leq s,t\leq k-1, t>k-s-2\\
\bbR,& \mbox{otherwise}
\end{array}
\right.
\end{equation}
(iii) If $\mu-\lambda-\nu$ is not a positive integer, then
\begin{equation} \nonumber \mathrm
H^1(\Sl(2,\bbR);  \cD_{\lambda,\nu;\mu})=0.
\end{equation}
\end{thm}
To proof the theorem above, we are required to proof the following
two Lemmas.
\begin{lem}
\label{lem1} Let $Y$ be a vector fields in $\Sl(2,\bbR)$ and let
$c:\cF_{\lambda}\otimes \cF_{\nu}\rightarrow \cF_{\mu}$ be a
bilinear differential operator defined as follows. For all
$\phi\in\cF_{\lambda}$ and for all $\psi\in\cF_{\nu}:$
$$
c(\phi,\psi)=\sum_{i+j=k} \alpha_{i,j}\,Y'\,\phi^{(i)}\,\psi^{(j)}
+\sum_{i+j=k-1}\beta_{i,j}\,Y''\,\phi^{(i)}\,\psi^{(j)}, $$ where
$\alpha_{i,j}$ and $\beta_{i,j}$ are constants. Then, for all
$X\in \Sl(2,\bbR),$ we have
$$
\begin{array}{ccl}
L_X^{\lambda,\nu;\mu} c(\phi,\psi)&\!\!\!\!=\!\!\!\!&\displaystyle
-\frac{1}{2}\,Y'X'' \sum_{i+j=k-1}\left((i+1)\left (i+2\lambda
\right )\alpha_{i+1,j}+(j+1)\left (j+2\nu
\right) \alpha_{i,j+1}\right )\,\phi^{(i)}\psi^{(j)}\\[2mm]
&& + \displaystyle X'Y''\sum_{i+j=k-1}(\mu-\lambda-\nu -i-j)\,
\beta_{i,j}\,\phi^{(i)}\,\psi^{(j)}\displaystyle + XY''
\sum_{i+j=k} \alpha_{i,j}\,\phi^{(i)}\psi^{(j)}
\end{array}
$$
\end{lem}
{\bf Proof.} Straightforward computation using the definition
(\ref{act}).
\begin{lem}
\label{lem2} Let $c:\cF_{\lambda}\otimes \cF_{\nu}\rightarrow
\cF_{\mu}$ be a bilinear differential operator defined as follows.
For all $\phi\in\cF_{\lambda}$ and for all $\psi\in\cF_{\nu}:$
$$
c(\phi,\psi)=\sum_{i+j=k} c_{i,j}\,\phi^{(i)}\,\psi^{(j)},
$$ where $c_{i,j}$ are constants. Then, for all $X\in \Sl(2,\bbR),$ we have
$$
L_X^{\lambda,\nu;\mu} c(\phi, \psi)=\frac{1}{2}\sum_{i+j=k-1}\left
( (i+1)(2i+\lambda)\,c_{i+1,j}+ (j+1)(2j+\nu)\,c_{i,j+1}\right
)X'' \phi^{(i)}\,\psi^{(j)}
$$
\end{lem}
{\bf Proof.} Straightforward computation using the definition
(\ref{act}).\\

Now we are in position to prove Theorem (\ref{leco}). Any
1-cocycle on $\Sl(2,\bbR)$ should retains the following general
form:
\begin{equation}
\label{sam}
c(X,\phi,\psi)=\sum_{j+j=k}\alpha_{i,j}\,X'\phi^{(i)}\psi^{(j)}+
\sum_{j+j=k-1}\beta_{i,j}\,X''\phi^{(i)}\psi^{(j)},
\end{equation}
where $\alpha_{i,j}$ and $\beta_{i,j}$ are constants. The higher
degree terms on $X$ are absent from the formula above, as they
vanish on the Lie algebra $\Sl(2,\bbR).$

The 1-cocycle condition reads as follows: for all $\phi \in
\cF_{\lambda},$ for all $\psi\in \cF_{\nu}$ and for all
$X\in\Sl(2,\bbR),$ we have
$$
c([X,Y],\phi,\psi)
-L_X^{\lambda,\nu;\mu}\,c(Y,\phi,\psi)+L_Y^{\lambda,\nu;\mu}\,
c(X,\phi,\psi)=0.
$$
A direct computation, and by using Lemma \ref{lem1}, proves that
the coefficient of the component $\phi^{(m)}\,\psi^{(n)}$ in the
1-cocycle condition above is equal to
\begin{equation}
\label{sl} \frac{1}{2}\,(Y'X''-X'Y'') \left ((m+1)\left(m
+2\lambda \right )\alpha_{m+1,n}+(n+1)\left(n+2\nu \right
)\alpha_{m,n+1}\right ).
\end{equation}
The formula (\ref{sl}) turns into zero once restricted to the
affine Lie algebra $\mbox{Span}\{\frac{d}{dx},x\frac{d}{dx}\}.$ We
are required, therefore, to study the annihilation of the formula
(\ref{sl}) for the two vector fields $X=x\frac{d}{dx}$ and
$Y=x^2\frac{d}{dx}.$ For these vector fields, the 1-cocycle
property will be equivalent to the system
\begin{equation}
\label{lam} (m+1)\left(m+2\lambda \right
)\alpha_{m+1,n}+(n+1)\left(n+2\nu \right )\alpha_{m,n+1}=0,
\end{equation}
where $m+n=k-1.$

Now we are going to deal with trivial 1-cocycles, and show how the
general 1-cocycles (\ref{sam}) can be eventually trivial. Any
trivial 1-cocycle should be of the form
$$
L_X^{\lambda,\nu;\mu}c,
$$
where $c$ is an operator $c:\cF_{\lambda}\otimes
\cF_{\nu}\rightarrow \cF_{\mu}$ defined as $c(\phi,
\psi)=\sum_{i+j=k}c_{i,j}\phi^{(i)}\psi^{(j)}.$ By using Lemma
\ref{lem2}, we have
\begin{equation}
\label{tya} L_x^{\lambda,\nu;\mu}c=\frac{1}{2}\,X''\sum_{m+n=k-1}
\left (-(m+1)\left(m +2\lambda \right )c_{m+1,n}-(n+1)\left(n
+2\nu \right )c_{m,n+1}\right )\phi^{(m)}\,\psi^{(n)}.
\end{equation}
We emphasize on the fact that the component $X'$ is absent from
the formula above. To complete the proof we distinguish many
cases:

(i) If $\lambda\not=-\frac{s}{2}$ and $\nu\not=-\frac{t}{2},$
where $s,t\in \{0,\ldots,k-1\},$ then the space of solutions of
the system (\ref{lam}) is one-dimensional; it is generated by
$\alpha_{0,k}.$ Now, we will explain how the constant
$\beta_{i,j}$ can be eliminated from our initial cocycle
(\ref{sam}). We add the coboundary $L_X^{\lambda,\nu;\mu}c $ of
equation (\ref{tya}) to our 1-cocycle (\ref{sam}). The constants
$c_{i,j}$ are chosen such that
$$
\beta_{m,n}=-2((m+1)\left(m+2\lambda \right
)c_{m+1,n}+(n+1)\left(n +2\nu) \right )c_{m,n+1}.
$$
This requirement is always possible, as $\lambda\not=-\frac{s}{2}$
and $\nu\not=-\frac{t}{2}.$ Therefore, our 1-cocycle (\ref{sam})
should only contain components in $X'.$ Now, by using (\ref{lam})
we can see that the cohomology group in question is
one-dimensional, generated by the following 1-cocycle.
\begin{eqnarray}
{\mathfrak a}(X,\phi,\psi)&=&X'\phi\,\psi^{(k)}\\
&&+\sum_{u+v=k-1}
(-1)^{k-v}\binom{k}{v}\frac{(v+2\nu)(v+1+2\nu)\cdots (k-1+2\nu)
}{(u+2\lambda)(u-1+2\lambda)\cdots
(2\lambda)}\,X'\phi^{(u+1)}\psi^{(v)}. \nonumber
\end{eqnarray}
(ii) If $\nu=-\frac{t}{2}$ and $\lambda\not=-\frac{s}{2},$ then
the constants
$\alpha_{k-t,t},\alpha_{k-t+1,t-1},\ldots,\alpha_{k,0}$ are zero
and the space of solutions of the system (\ref{lam}) is
one-dimensional, generated by $\alpha_{0,k}.$ The constant
$\beta_{i,j}$ can be eliminated by the same method as in Part (i).
We have just proved that the cohomology group in question is
one-dimensional, generated by the 1-cocycle:
\begin{equation}
{\mathfrak
b}(X,\phi,\psi)=X'\phi\,\psi^{k}+\sum_{u+1+v=k}\alpha_{u+1,v}
\,X'\phi^{(u+1)}\psi^{(v)},
\end{equation}
where
$$
\alpha_{u+1,v}=\left\{
\begin{array}{ll}
0,& \mbox{if} \quad v\leq t\\
(-1)^{k-v}\binom{k}{v}\frac{(v+2\nu)(v+1+2\nu)\cdots (k-1+2\nu)
}{(u+2\lambda)(u-1+2\lambda)\cdots (2\lambda)},&\mbox{otherwise}
\end{array}
\right.
$$
(iii) If $\lambda=-\frac{s}{2}$ and $\nu\not=-\frac{t}{2},$ then -
and as in Part (ii) - the cohomology group in question is
one-dimensional, generated by the 1-cocycle:
\begin{equation}
{\mathfrak
c}(X,\phi,\psi)=X'\phi^{k}\,\psi+\sum_{u+1+v=k}\alpha_{u,v+1}
\,X'\phi^{(u)}\psi^{(v+1)},
\end{equation}
where
$$
\alpha_{u,v+1}=\left\{
\begin{array}{ll}
0,& \mbox{if} \quad u\leq s\\
(-1)^{k-u}\binom{k}{u}\frac{(u+2\lambda)(u+1+2\lambda)\cdots
(k-1+2\lambda) }{(v+2\nu)(v-1+2\nu)\cdots
(2\nu)},&\mbox{otherwise}
\end{array}
\right.
$$
(iv) If $\lambda=-\frac{s}{2}$ and $\nu=-\frac{k-s-1}{2},$ where
$s\in \{0,\ldots,k-1\},$ then the space of solutions of the system
(\ref{lam}) is two dimensional; it is generated by
$\alpha_{s+1,k-s-1}$ and $\alpha_{s,k-s}.$ Now, we will explain
how the constant $\beta_{i,j}$ (but not $\beta_{s,k-s-1}$) can be
eliminated. We add the coboundary $L_X^{\lambda,\nu}c $ of
(\ref{tya}) to our 1-cocycle (\ref{sam}). The constant $c_{i,j}$
are chosen such that
$$
\beta_{m,n}=-2((m+1)\left(m+2\lambda \right
)c_{m+1,n}+(n+1)\left(n+2\nu )\right )c_{m,n+1}.
$$
This requirement is always satisfied, except for $\beta_{s,k-s-1}$
because the component $\phi^{s}\,\psi^{k-s-1}$ of our trivial
1-cocycle (\ref{tya}) has a trivial coefficient. Finally, we have
just proved that the cohomology group in question is
three-dimensional, generated by the 1-cocycles:
\begin{equation}
\begin{array}{ccl}
{\mathfrak d}(X,\phi,\psi)&=&\beta_{s,k-s-1}\,
X''\phi^{s}\,\psi^{k-s-1}+\alpha_{0,k}\,X'\phi\,\psi^{k}+
\alpha_{k,0}\,X'\phi^{k}\,\psi\\[2mm]
&&\displaystyle +\sum_{\substack{u+v=k\\u,v\not= 0}}\alpha_{u,v}
\,X'\phi^{(u)}\psi^{(v)},
\end{array}
\end{equation}
where
$$
\alpha_{u,v}=\left\{
\begin{array}{ll}
\displaystyle(-1)^{k-v}\binom{k}{v}\frac{(v+2\nu)(v+1+2\nu)\cdots
(k-1+2\nu) }{(u-1+2\lambda)(u-2+2\lambda)\cdots
(2\lambda)}\alpha_{0,k},& \mbox{if}
\quad u\leq s\\[3mm]
\displaystyle
(-1)^{k-u}\binom{k}{u}\frac{(u+2\lambda)(u+1+2\lambda)\cdots
(k-1+2\lambda) }{(v-1+2\nu)(v-2+2\nu)\cdots
(2\nu)}\alpha_{k,0},&\mbox{if}\quad u\geq s+1
\end{array}
\right.
$$
(v) If $\lambda=-\frac{s}{2}$ and $\nu=-\frac{t}{2},$ where
$s,t\in \{0,\ldots,k-1\}$ but $ t\leq k-s-2,$ then the space of
solutions of the system (\ref{lam}) is one-dimensional, generated
by $\alpha_{s+1,k-s-1}.$ The constant $\beta_{i,j}$ can be
eliminated as explained before. Thus, the cohomology group in
question is one-dimensional, generated by the 1-cocycle:
\begin{equation}
\begin{array}{ccl}
{\mathfrak e}
(X,\phi,\psi)=X'\,\phi^{(s+1)}\,\psi^{(k-s-1)}\displaystyle
+\sum_{\substack{u+v=k\\u\not= s+1}}\alpha_{u,v}
\,X'\phi^{(u)}\psi^{(v)},
\end{array}
\end{equation}
where
$$
\alpha_{u,v}=\left\{
\begin{array}{ll}
0,&\mbox{if}\quad u\leq s\\
0,&\mbox{if}\quad v\leq t\\
(-1)^{k-s-1-v}\frac{(v+1)(v+2)\cdots n }{u(u-1)\cdots s
}\frac{(v+2\nu)(v+1+2\nu)\cdots (k-s-2+2\nu)
}{(u-1+2\lambda)(u-2+2\lambda)\cdots (s+1+2\lambda)},&
\mbox{otherwise}
\end{array}
\right.
$$
(vi) If $\lambda=-\frac{s}{2}$ and $\nu=-\frac{t}{2},$ where
$s,t\in \{0,\ldots,k-1\}$ but $ t>k-s-2,$ then the space of
solutions of the system (\ref{lam}) is two-dimensional, generated
by $\alpha_{s+1,k-s-1}$ and $\alpha_{k-t-1,t+1}.$ The constant
$\beta_{i,j}$ can be eliminated as explained before, except for
$\beta_{k-t-1,t}.$ Thus, the cohomology group in question is
three-dimensional, generated by the 1-cocycle:
\begin{equation}
\begin{array}{ccl}
{\mathfrak f}(X,\phi,\psi)&=&\beta_{k-t-1,t}\,
X''\phi^{k-t-1}\,\psi^{t}+\alpha_{0,k}\,X'\phi\,\psi^{k}+
\alpha_{k,0}\,X'\phi^{k}\,\psi\\[2mm]
&&\displaystyle +\sum_{\substack{u+v=k\\u,v\not= 0}}\alpha_{u,v}
\,X'\phi^{(u)}\psi^{(v)},
\end{array}
\end{equation}
where
$$
\alpha_{u,v}=\left\{
\begin{array}{ll}
\displaystyle(-1)^{k-v}\binom{k}{v}\frac{(v+2\nu)(v+1+2\nu)\cdots
(k-1+2\nu) }{(u-1+2\lambda)(u-2+2\lambda)\cdots
(2\lambda)}\alpha_{0,k},& \mbox{if}
\quad v\geq t+1\\[3mm]
\displaystyle
(-1)^{k-u}\binom{k}{u}\frac{(u+2\lambda)(u+1+2\lambda)\cdots
(k-1+2\lambda) }{(v-1+2\nu)(v-2+2\nu)\cdots
(2\nu)}\alpha_{k,0},&\mbox{if}\quad u\geq s+1
\end{array}
\right.
$$
\begin{rmk}{\rm
The first cohomology group of the Lie algebra $\Sl(2,\mathbb{R})$
with coefficients in the space of linear differential operators
has been computed in \cite{l}. The explicit 1-cocycles that span
this cohomology group has first arisen in \cite{gar}.}
\end{rmk}
\section{$\Sl(2,\bbR)$-invariant differential operators}
In this section we will investigate differential operators on
tensor densities that are $\Sl(2,\mathbb{R})$-invariant. These
results will be useful for the computation of cohomology.
\begin{pro}\cite{g}
\label{gor}
There exist unique (up to constants)
$\Sl(2,\mathbb{R})$-invariant bilinear differential operators
$J_k^{\l,\mu}: \cF_\l\otimes \cF_\nu \rightarrow \cF_{\l+\nu+k}$
given by
\begin{equation}
\label{hmida} J_k^{\lambda,\nu}(\phi, \psi ) =\sum_{i+j=k}
c_{i,j}\,\phi^{(i)}\,\psi^{(j)},
\end{equation}
where the constants $c_{i,j}$ are characterized as follows:

(i) If $\lambda, \nu \not\in \{0,-1/2,\ldots,-s/2,\ldots\},$ the
coefficients $c_{i,j}$ are given by
$$
c_{i,j}=\displaystyle (-1)^{i} {k \choose
i}\frac{(2\l-i)(2\lambda -i-1)\ldots (2\lambda
-k+1)}{(2\nu-j)(2\nu-j+1)\ldots (2\nu -k+1)}.
$$

(ii) If $\lambda$ or $\nu \in \{0,-1,-1/2,\ldots, -s/2,\ldots\},$
the coefficients $c_{i,j}$ satisfy the recurrence relation
\begin{equation}
\label{ak}
(i+1)(i+2\lambda)\,c_{i+1,j}+(j+1)(j+2\nu)\,c_{j,i+1}=0.
\end{equation}
Moreover, the space of solutions of the system (\ref{ak}) is
two-dimensional if $\lambda=-\frac{s}{2}$ and $\nu=-\frac{t}{2}$
but $t>k-s-2,$ and one-dimensional otherwise.
\end{pro}

\begin{pro}
\label{inv} There exist $\Sl(2,\mathbb{R})$-invariant trilinear
differential operators $K_k^{\l,\nu,\tau}: \cF_\l\otimes
\cF_\nu\otimes \cF_{\tau} \rightarrow \cF_{\l+\nu+\tau+k}$ given
by
\begin{equation}
\label{hmida2} K_k^{\lambda,\nu, \tau}(\phi, \varphi, \psi )
=\sum_{i+j+l=k} c_{i,j,l}\,\phi^{(i)}\,\varphi^{(j)}\,\psi^{(l)} ,
\end{equation}
where the constants $c_{i,j,l}$ are characterized by the
recurrence formula
\begin{equation}
\label{ars}
i(i-1+2\lambda)\,c_{i,j,l}+(j+1)(j+2\nu)\,c_{i-1,j+1,l}+
(l+1)(l+2\tau)\,c_{i-1,j,l+1}=0,
\end{equation}
where $i+j+l=k.$

If $\lambda, \nu$ and $\tau$ are generic, then the space of
solutions is $(k+1)$-dimensional.
\end{pro}
\begin{pro}
\label{inv2} There exist $\Sl(2,\mathbb{R})$-invariant trilinear
differential operators\\ $K_k^{\nu,\tau}: \Vect(\mathbb{RP}^1)
\otimes \cF_\nu \otimes \cF_{\tau} \rightarrow \cF_{\nu+\tau+k-1}$
that vanishe on $\Sl(2,\mathbb{R})$ given by
\begin{equation}
\label{hmida2} K_k^{\nu, \tau}(X, \varphi, \psi ) =\sum_{i+j+l=k}
c_{i,j,l}\,X^{(i)}\,\varphi^{(j)}\,\psi^{(l)} ,
\end{equation}
where the constants $c_{i,j,l}$ are as in (\ref{ars}) but
$c_{0,j,k-j}=c_{1,j,k-j-1}=c_{2,j,k-j-2}=0.$ Moreover, the space
of solutions is $(k-2)$-dimensional, for all $\nu$ and $\tau.$
\end{pro}
{\bf Proof of Proposition (\ref{inv}) and (\ref{inv2}).} We are
going to prove Proposition (\ref{inv}) and (\ref{inv2})
simultaneously. Any differential operator $K_k^{\l,\nu,\tau}:
\cF_\l\otimes \cF_\nu\otimes \cF_{\tau} \rightarrow \cF_{\mu}$ is
of the form
\begin{equation}
\nonumber K_k^{\lambda,\nu, \tau}(\phi, \varphi, \psi )
=\sum_{i+j+l=k} c_{i,j,l}\,\phi^{(i)}\,\varphi^{(j)}\,\psi^{(l)} ,
\end{equation}
where $c_{i,j,l}$ are functions.

The $\Sl(2,\mathbb{R})$-invariant property of the operators
$K_k^{\l,\nu,\tau}$ reads as follows.
$$
L_X^\nu K_k^{\lambda,\nu, \tau}(\phi, \varphi, \psi
)=K_k^{\lambda,\nu, \tau}(L_X^\lambda\phi, \varphi, \psi )+
K_k^{\lambda,\nu, \tau}(\phi, L_X^\nu\varphi, \psi
)+K_k^{\lambda,\nu, \tau}(\phi, \varphi, L_X^\tau\psi ).
$$
The invariant property with respect to the affine Lie algebra
$\mathrm{Span}\{\frac{d}{dx}, x\frac{d}{dx}\}$ implies that
$c'_{i,j,l}=0$ and $\mu=\l+\nu+\tau+k.$ On the other hand, the
invariant property with respect to the vector fields
$x^2\frac{d}{dx}$ is equivalent to the system (\ref{ars}). If
$\lambda,$ $\nu$ and $\tau$ are generic, then the space of
solutions is $(k+2)$-dimensional, generated by $c_{k-1,1,0},
c_{k-1,0,1}, c_{k-2,2,0}, c_{k-3,3,0},\ldots, c_{0,k,0}.$

Now, the proof of Proposition (\ref{inv2}) follows as above by
putting $\lambda=-1.$ In this case, the space of solutions is
$(k-2)$-dimensional, spanned by $c_{3,k-3,0}, c_{3,k-4,1},
c_{3,k-5,2},\ldots, c_{3,0,k-3}.$
\section{Cohomology of $\Vect(\mathbb{RP}^1)$ acting on
${\cal D}_{\lambda,\nu;\mu}$}
In this section, we will compute the first cohomology group of
$\Vect(\mathbb{RP}^1)$ with values in $\cD_{\lambda,\nu;\mu},$
vanishing on $\Sl(2,\mathbb{R}).$
\begin{thm}
\label{main} (i) If $\mu-\lambda-\nu=2,$ then
\begin{equation} \nonumber
\mathrm H^1(\Vect(\mathbb{RP}^1),\Sl(2,\bbR);
\cD_{\lambda,\nu;\mu})= \left\{
\begin{array}{ll}
\bbR,& \mbox{if} \quad \lambda+\nu+1=0\\
\bbR,& \mbox{if} \quad \lambda=0 \mbox{ together with }
\nu\not =-\frac{1}{2}, \mbox{ and vice versa}\\
0,& \mbox{otherwise}
\end{array}
\right.
\end{equation}
(ii) If $\mu-\lambda-\nu=3,$ then
\begin{equation} \nonumber \mathrm H^1(\Vect(\mathbb{RP}^1),\Sl(2,\bbR);
\cD_{\lambda,\nu;\mu})= \left\{
\begin{array}{lll}
\bbR^2,& \mbox{ if }\quad  (\lambda,\nu) = &(0,0), (-2,0), (0, -2),
(-\frac{2}{3}, -\frac{2}{3})\\[2mm]
\bbR,&\mbox{ if }\quad  (\lambda,\nu)\not =
&(-\frac{1}{2},-1),(-1,-\frac{1}{2}), (-1,-1),\\[2mm]
& \mbox{  }\quad &  (0,-1), (-1,0),
(-\frac{1}{2},-\frac{1}{2})\\
0,& \mbox{otherwise}
\end{array}
\right.
\end{equation}
(iii) If $\mu-\lambda-\nu=4,$ then
\begin{equation} \nonumber \mathrm H^1(\Vect(\mathbb{RP}^1),\Sl(2,\bbR);
\cD_{\lambda,\nu;\mu})= \left\{
\begin{array}{ll}
\bbR^2,& \mbox{ if }\quad  (\lambda,\nu) \not= (0,-\frac{3}{2}),
(-\frac{3}{2},0),(-1,-1)\\[2mm]
& (-\frac{1}{2}, -\frac{3}{2}), (-\frac{3}{2}, -\frac{1}{2})
,(-\frac{3}{2},-1),(-1, -\frac{3}{2}),\\[2mm]
&(-\frac{3}{2}, -\frac{3}{2}), (-\frac{1}{2},-1), (-1,
-\frac{1}{2})\\[2mm]
\bbR,& \mbox{otherwise}
\end{array}
\right.
\end{equation}
(iv) If $\mu-\lambda-\nu=5,$ then
\begin{equation} \nonumber
\mathrm H^1(\Vect(\mathbb{RP}^1),\Sl(2,\bbR);
\cD_{\lambda,\nu;\mu})= \left\{
\begin{array}{ll}
\bbR^3,& \mbox{if} \quad (\lambda,\nu) = (0,0),
(-4,0),(0,-4)\\[2mm]
\bbR^2,& \mbox{if} \quad (\lambda,\nu)\not =
(-\frac{1}{2},-\frac{3}{2}),(-\frac{3}{2},-\frac{1}{2}),
(-\frac{1}{2},-2),\\[2mm]& (-2,-\frac{1}{2}),
(-1,-2) (-2,-1),(-1,-\frac{3}{2}),\\[2mm]& (-\frac{3}{2},
-1),(-1,-1)(-\frac{3}{2},-2), (-2,-\frac{3}{2}),\\[2mm] &
(-\frac{3}{2},-\frac{3}{2})\\[2mm]
 \bbR,& \mbox{otherwise}
\end{array}
\right.
\end{equation}
(v) If $\mu-\lambda-\nu=6,$ then
\begin{equation} \nonumber
\mathrm H^1(\Vect(\mathbb{RP}^1),\Sl(2,\bbR);
\cD_{\lambda,\nu;\mu})= \left\{
\begin{array}{ll}
\bbR^3,& \mbox{if} \quad (\lambda,\nu) = (
\frac{-5\pm\sqrt{19}}{2},0), (0,
\frac{-5\pm\sqrt{19}}{2}),\\[2mm]
&(\frac{\sqrt{19}-5}{2}, \frac{-\sqrt{19}-5}{2}),
(\frac{-\sqrt{19}-5}{2},\frac{\sqrt{19}-5}{2})\\[2mm]
\bbR^2,& \mbox{if} \quad (\lambda,\nu) \not= (-2,-\frac{1}{2}),
(-\frac{1}{2},-2),
(-\frac{5}{2},-\frac{1}{2}),\\[2mm]
&(-\frac{1}{2},-\frac{5}{2}),(-1,-\frac{3}{2}),
(-\frac{3}{2},-1),(-1,-2),\\[2mm]
& (-2,-1),(-1,-\frac{5}{2}),(-\frac{5}{2},-1),
(-\frac{3}{2},-\frac{3}{2}),\\[2mm]
& (-\frac{3}{2},-2),(-2,-\frac{3}{2}),(-\frac{3}{2},-\frac{5}{2}),
(-\frac{5}{2},\frac{3}{2}),\\[2mm]
& (-2,-2),
(-2, -\frac{5}{2}), (-\frac{5}{2}, -2),\\[2mm]
 \bbR,& \mbox{otherwise}
\end{array}
\right.
\end{equation}
(vi) If $\mu-\lambda-\nu=7,$ then
\begin{equation} \nonumber
\mathrm H^1(\Vect(\mathbb{RP}^1),\Sl(2,\bbR);
\cD_{\lambda,\nu;\mu})= \left\{
\begin{array}{ll}
\bbR^2,& \mbox{if} \quad (\lambda,\nu) =
(\frac{-5-\sqrt{19}}{2},\frac{-5-\sqrt{19}}{2}),
(\frac{\sqrt{19}-5}{2},\frac{\sqrt{19}-5}{2}),\\[2mm]
&(0,\nu), (\lambda,0), (\lambda, -6-\lambda),
(\sqrt{19}-1,\frac{-5-\sqrt{19}}{2}),\\[2mm]
&(\frac{-5-\sqrt{19}}{2},\sqrt{19}-1),
(-\sqrt{19}-1,\frac{-5+\sqrt{19}}{2}),\\[2mm]
&(\frac{-5+\sqrt{19}}{2},-\sqrt{19}-1)\\[2mm]
\bbR,& \mbox{if} \quad (\lambda,\nu) \not =  (-\frac{5}{2},
-\frac{1}{2}),
 (-\frac{1}{2}, -\frac{5}{2}), (-\frac{1}{2}, -3),
 \\[2mm]
& (-3,-\frac{1}{2}),(-1,-2), (-2,-1),(-1,-3),\\[2mm]
 &  (-3,-1), (-\frac{3}{2}, -2),(-2, -\frac{3}{2}),
 (-2,-1),\\[2mm]
 &(-1,-2),(-2,-\frac{3}{2}), (-\frac{3}{2}, -2),
 (-\frac{3}{2},-\frac{5}{2}), \\[2mm]
 &(-2,-3),(-3,-2),(-\frac{5}{2},-\frac{3}{2}),
 (-2,-\frac{5}{2}),\\[2mm]
& (-\frac{5}{2}, -2), (-3,-\frac{5}{2}),
(-\frac{5}{2}, -3) \\[2mm]
 0, & \mbox{otherwise}
\end{array}
\right.
\end{equation}
(vii) If $\mu-\lambda-\nu$ is not like above but $\lambda$ and
$\nu$ are generic then
\begin{equation} \nonumber
\mathrm H^1(\Vect(\mathbb{RP}^1),\Sl(2,\bbR);
\cD_{\lambda,\nu;\mu})=0.
\end{equation}
\end{thm}
\section{Proof of Theorem (\ref{main})}
To proof Theorem (\ref{main}) we proceed bye following the three
steps:
\begin{enumerate}
\item We will investigate the dimension of the space of operators
that satisfy the 1-cocycle condition. By Proposition (\ref{inv2}),
its dimension is at most $k-2,$ where $k=\mu-\lambda-\nu,$ since
any 1-cocycle that vanishes on $\Sl(2,\mathbb{R})$ is certainly
$\Sl(2,\mathbb{R})$-invariant (cf. \cite{bo1,lo2}). \item We will
study all trivial 1-cocycles, namely, operators of the form
$$
L_XB,
$$
where $B$ is a bilinear operator. As our 1-cocycles vanish on the
Lie algebra $\Sl(2,\mathbb{R}),$ it follows that the operator $B$
coincides with the transvectant $J_k^{\lambda,\nu}.$ By using
Proposition (\ref{gor}), we will determine different values of
$\lambda$ and $\nu$ for which the space of operators of the form
$L_X J_k^{\lambda,\nu}$ is zero, one or two-dimensional.

\item By taking into account Part 1 and part 2 - and depending on
$\lambda$ and $\nu$ - the dimension of the cohomology group
$\mathrm{H}^1(\Vect(\mathbb{RP}), \Sl(2,\mathbb{R});
\cD_{\lambda,\nu;\mu})$ will be equal to
$$
\mbox{dim(operators that are 1-cocycles)}-\mbox{dim(operators of
the form } L_XJ_{k}^{\lambda,\nu})
$$
\end{enumerate}
We need also the following Lemma.
\begin{lem}
\label{fun} Every 1-cocycle on $\mathrm{Vect}(\mathbb{RP}^1)$ with
values in $\cD_{\lambda,\nu;\mu}$ is differentiable.
\end{lem}
{\bf Proof.} See \cite{lo2}.

Now we are in position to prove Theorem (\ref{main}). By Lemma
(\ref{fun}), any 1-cocycle on $\Vect(\mathbb{RP}^1)$ should
retains the following general form:
\begin{equation}
\label{sam2}
c(X,\phi,\psi)=\sum_{j+j+l=k}\alpha_{l,i,j}\,X^{(l)}\phi^{(i)}\psi^{(j)},
\end{equation}
where $\alpha_{l,i,j}$ are constants. The fact that this 1-cocycle
vanishes on $\Sl(2,\mathbb{R})$ implies that
$$
c_{0,i,j}=c_{1,i,j}=c_{2,i,j}=0.
$$ The 1-cocycle condition reads
as follows: for all $\phi \in \cF_{\lambda},$ for all $\psi\in
\cF_{\nu}$ and for all $X\in\Vect(\mathbb{RP}^1),$ one has
$$
c([X,Y],\phi,\psi)
-L_X^{\lambda,\nu;\mu}\,c(Y,\phi,\psi)+L_Y^{\lambda,\nu;\mu}\,
c(X,\phi,\psi)=0.
$$
\subsection{The case when $k=2$}
If $\mu-\lambda-\nu=2,$ equation (\ref{sam2}) shows that only one
1-cocycle spans the cohomology group of Theorem (\ref{main}); it
is given by
\begin{equation}
\label{co3}
{\mathfrak L}(X,\phi,\psi):=X'''\,\phi\,\psi.
\end{equation}
Let us study the triviality of this 1-cocycle. A direct
computation proves that
$$
L_XJ_2^{\lambda,\nu}(\phi,\psi)=(-\lambda c_{2,0}-\nu
c_{0,2})X'''\,\phi\,\psi,
$$
\begin{enumerate}
\item If $\lambda=0$ and $\nu\not= -\frac{1}{2},$ then by Part (ii) of
Proposition (\ref{gor}), $c_{0,2}=0$ and the 1-cocycle is not
trivial. The result holds when $\nu=0$ and $\lambda\not=
-\frac{1}{2}$ as well.
\item If $\lambda=0$ and $\nu=-\frac{1}{2},$ then
$c_{0,2}=-\frac{1}{\nu}$ and therefore
$L_XJ_2^{0,-\frac{1}{2}}(\phi, \psi)={\mathfrak L}(X,\phi, \psi).$
Hence, the 1-cocycle is trivial. The result holds true when
$\nu=0$ and $\lambda= -\frac{1}{2}.$
\item If $\lambda$ and $\nu$ are not like above, then Part (i) of
Proposition (\ref{gor}) implies that
$$
-\lambda \,c_{2,0}-\nu\,
c_{0,2}=-2\nu\frac{\lambda+\nu+1}{1+2\lambda}c_{1,1},
$$
where $c_{1,1}\not =0.$ Thus, for $\lambda+\nu+1=0$ the 1-cocycle
(\ref{co3}) is not trivial; otherwise, the 1-cocycle is trivial.
\end{enumerate}
\subsection{The case when $k=3$}
If $\mu-\lambda-\nu=3,$ equation (\ref{sam2}) shows that the
1-cocycles that span the cohomology group of Theorem (\ref{main})
are of the form
\begin{equation}
\label{co3} {\mathfrak
M}(X,\phi,\psi):=\gamma\,X^{(4)}\,\phi\,\psi+(\alpha_{1}\,\phi'\,\psi
+\alpha_{2}\,\phi\,\psi')\,X''',
\end{equation}
where $\gamma, \alpha_1$ and $\alpha_2$ are constants. The
1-cocycle condition implies that
$$
\gamma=-\frac{\lambda\,\alpha_1+\nu\,\alpha_2}{2}.
$$
 Let us study the triviality
of this 1-cocycle. A direct computation proves that
$$
\begin{array}{lcl}
L_XJ_3^{\lambda,\nu}(\phi,\psi)&=&(-\lambda \,c_{3,0}-\nu\,
c_{0,3})X^{(4)}\,\phi\,\psi-((1+3\lambda)\,c_{3,0}+\nu\,c_{1,2})
X'''\phi'\,\psi\\[2mm]
&&-((1+3\nu)\,c_{0,3}+\lambda\,c_{2,1}) X'''\phi\,\psi'
\end{array}
$$
\begin{enumerate}
\item Here we will characterize all values of $\lambda$ and $\nu$
for which $L_X J_3^{\lambda,\nu}=0.$ An easy computation using
Proposition (\ref{gor}) proves that these values are $(0,0),
(-2,0), (0,-2),$ and $(-\frac{2}{3},-\frac{2}{3}).$  If
$(\lambda,\nu)=(0,0),$ then $c_{0,3}=c_{3,0}=0$ and
$c_{1,2}=-c_{2,1}.$  It follows that the cohomology group is
two-dimensional, spanned by two 1-cocycles given as in (\ref{co3})
for $(\alpha_1,\alpha_2)=(1,0)$ and $(0,1).$ The results holds
true for $(\lambda,\nu)=(-2,0), (0,-2),$ and
$(-\frac{2}{3},-\frac{2}{3}).$

\item Here we will characterize all values of $\lambda$ and $\nu$
for which $L_X J_3^{\lambda,\nu}$ is generated by two parameters.
An easy computation using Proposition (\ref{gor}) proves that
these values are $(-\frac{1}{2},-1), (-1,-\frac{1}{2}),
(0,-1),(-1,0),(-\frac{1}{2},-\frac{1}{2}),$ and $(-1,-1).$ If
$(\lambda,\nu)=(-\frac{1}{2},-1),$ then $c_{1,2}=0,$ and
$c_{3,0}=\frac{2}{3}c_{2,1};$ however, $c_{0,3}$ and $c_{2,1}$ are
arbitrary. Therefore, the constants $c_{0,3}$ and $c_{2,1}$ can be
chosen such that
$$
L_X J^{\lambda,\nu}_3(\phi,\psi)= {\mathfrak M}(X,\phi,\psi).
$$
It follows that the cohomology group is trivial. The result holds
true for $(0,-1),(-1,0),$
$(-\frac{1}{2},-\frac{1}{2}),(-1,-\frac{1}{2}),$ and $(-1,-1).$
\item If $\lambda$ and $\nu$ are not like above, then the
transvectant $J_3$ is unique by Proposition (\ref{gor}). Whatever
the weights $\lambda$ and $\nu$ can take, the trivial 1-cocycle $
L_X J^{\lambda,\nu}_3$ is never identically zero. Therefore, one
of the constants $\gamma,$ $\alpha_1$ or $\alpha_2$ can be
eliminated by just adding the trivial 1-cocycle $ L_X
J^{\lambda,\nu}_3.$ Hence, the cohomology group is
one-dimensional.
\end{enumerate}
\subsection{The case when $k=4$}
If $\mu-\lambda-\nu=4,$ equation (\ref{sam2}) shows that the
1-cocycles that span the cohomology group of Theorem (\ref{main})
are of the form
\begin{equation}
\label{co4} {\mathfrak
N}(X,\phi,\psi):=\gamma\,X^{(5)}\,\phi\,\psi+X^{(4)}\,\left (
\alpha_{1}\,\phi'\,\psi +\alpha_{2}\,\phi\,\psi'\right )+X'''\left
(\beta_1\,\phi''\,\psi+\beta_2\,\phi\,\psi''+\beta_3\,\phi'\,
\psi' \right ),
\end{equation}
where $\gamma, \alpha_1, \alpha_2, \beta_1, \beta_2$ and $\beta_3$
are constants. The 1-cocycle condition is equivalent to the
following system
$$
\begin{array}{ll}
5\gamma=-\lambda\,\alpha_1-\nu\,\alpha_2,&
2\alpha_1=-(1+2\lambda)\,\beta_1-\nu \,\beta_3,\\
2\alpha_2=-(1+2\nu)\,\beta_2-\lambda\, \beta_3.
\end{array}
$$
The space of solutions of the system above is three-dimensional.
Let us study the triviality of the 1-cocycle (\ref{co4}). A direct
computation proves that
$$
\begin{array}{lcl}
L_XJ_4^{\lambda,\nu}(\phi,\psi)&\!\!\!\!\!\!=\!\!\!\!\!\!&(-\lambda
\,c_{4,0}-\nu\,
c_{0,4})X^{(5)}\,\phi\,\psi-((1+4\lambda)\,c_{4,0}+\nu\,c_{1,3})
X^{(4)}\phi'\,\psi\\[2mm]
&&-((1+4\nu)\,c_{0,4}+\lambda\,c_{3,1}) X^{(4)}\,\phi\,\psi'-
((4+6\lambda)\,c_{4,0}+\nu\, c_{2,2})\,X'''\,\phi''\,\psi\\[2mm]
&&-(((1+3\lambda)\,c_{3,1}+(1+3\nu)\, c_{1,3})\,\phi'\,\psi'-
((4+6\nu)\,c_{0,4}+\lambda\, c_{2,2})\,\phi\,\psi'')\,X'''
\end{array}
$$
\begin{enumerate}
\item For all values of $\lambda$ and $\nu$ one can easily prove
that the equation $L_XJ_4^{\lambda, \nu}=0$ has no solutions.

\item Now we will characterize all values of $\lambda$ and $\nu$
for which $L_X J_4^{\lambda,\nu}$ is generated by two parameters.
An easy computation using Proposition (\ref{gor}) proves that
these values are $(0,-\frac{3}{2}), (-\frac{3}{2},0),
(-\frac{1}{2},-\frac{3}{2}), (-\frac{3}{2},-\frac{1}{2}),
(-\frac{1}{2},-1),(-1,-\frac{1}{2}), (-1,-1), (-1,-\frac{3}{2}),$
$(-\frac{3}{2}, -1),$ and $(-\frac{3}{2}, -\frac{3}{2}).$ If
$(\lambda,\nu)=(0,-\frac{3}{2}),$ then the constants $c_{1,3}$ and
$c_{0,4}$ are arbitrary. On the other hand,
$$c_{2,2}=\frac{3}{2}\,c_{1,3}, \quad c_{3,1}=c_{3,1},\quad
c_{4,0}=\frac{1}{4}\,c_{1,3}.$$ The constant $c_{1,3}$ and
$c_{0,4}$ can be chosen in such a way that once adding the trivial
1-cocycle $L_XJ_4^{\lambda,\nu}$ above to our 1-cocycle
(\ref{co4}), the constants $\beta_1$ and $\beta_2$ disappear
completely. Hence, the cohomology group is one-dimensional. The
result holds true for the other values of $(\lambda,\nu).$ \item
If $\lambda$ and $\nu$ are not like above, then the transvectant
$J_4$ is unique by Proposition (\ref{gor}). Whatever the weights
$\lambda$ and $\nu$ can take, the trivial 1-cocycle $ L_X
J^{\lambda,\nu}_4$ is never identically zero. Therefore, one of
the constants $\gamma,$ $\alpha_1,$ $\alpha_2,$ $\beta_1,$
$\beta_2$ or $\beta_3$ can be eliminated by just adding the
trivial 1-cocycle $ L_X J^{\lambda,\nu}_4.$ Hence, the cohomology
group is two-dimensional.
\end{enumerate}
\subsection{The case when $k=5$}
If $\mu-\lambda-\nu=5,$ equation (\ref{sam2}) shows that the
1-cocycles that span the cohomology group of Theorem (\ref{main})
are of the form
\begin{eqnarray}
{\mathfrak
O}(X,\phi,\psi)&:=&\gamma\,X^{(6)}\,\phi\,\psi+X^{(5)}\,\left (
\alpha_{1}\,\phi'\,\psi +\alpha_{2}\,\phi\,\psi'\right
)+X^{(4)}\left
(\beta_1\,\phi''\,\psi+\beta_2\,\phi\,\psi''+\beta_3\,\phi'\,
\psi' \right )\nonumber\\
&&X'''\left
(\eta_1\,\phi^{(3)}\,\psi+\eta_2\,\phi^{(2)}\,\psi'+\eta_3\,\phi'\,
\psi^{(2)}+ \eta_4\,\phi\,\psi^{(3)}\right )\label{co5}
\end{eqnarray}
The 1-cocycle condition is equivalent to the following system
$$
\begin{array}{ll}
5\gamma=-\lambda\,\beta_1-\nu\,\beta_2+\lambda \,\eta_1+\nu\,
\eta_4,
&9\gamma=-\lambda\,\alpha_1-\nu\,\alpha_2,\\
5\alpha_1=-(1+2\lambda)\,\beta_1-\nu \,\beta_3,&
5\alpha_2=-(1+2\nu)\,\beta_2-\lambda\, \beta_3,\\
2\beta_1=-(3+3\lambda)\,\eta_1-\nu\,\eta_2,&
2\beta_2=-\lambda\,\eta_3-(3+3\nu)\,\eta_4\\
2\beta_3=-(1+2\lambda)\,\eta_2-(1+2\nu)\,\eta_3.
\end{array}
$$
The space of solutions of the system above is four-dimensional for
$(\lambda,\nu)=(0,0),$ $(0,-2),$ $(-2,0), $     $(0,-4),$$
(-4,0),$ and $(-2,-2),$ and three-dimensional otherwise. According
to these values, let us study the triviality of the 1-cocycle
(\ref{co5}). A direct computation proves that
$$
\begin{array}{lcl}
L_XJ_5^{\lambda,\nu}(\phi,\psi)&\!\!\!\!\!\!\!=\!\!\!\!\!\!\!
&(-\lambda
\,c_{5,0}-\nu\,
c_{0,5})X^{(6)}\,\phi\,\psi-((1+5\lambda)\,c_{5,0}+\nu\,c_{1,4})
X^{(5)}\phi'\,\psi\\[2mm]
&&-((1+5\nu)\,c_{0,5}+\lambda\,c_{4,1}) X^{(5)}\,\phi\,\psi'-
X^{(4)}(((5+10\lambda)\,c_{5,0}-\nu\, c_{2,3})\,\phi''\,\psi\\[2mm]
&&-((1+3\lambda)\,c_{4,1}+(1+3\nu)\, c_{1,4})\,\phi'\,\psi'-
((5+10\nu)\,c_{0,5}+\lambda\,
c_{2,3})\,\phi\,\psi'')\\[2mm]
&&-X'''((\nu\,c_{3,2}+(10+12\lambda)\,c_{5,0})\,\phi^{(3)}\,\psi
-(\nu\, c_{2,3}+(10+12\nu)\,c_{0,5})\,\phi\,\psi^{(3)}\\[2mm]
&&-((1+3\lambda)\, c_{3,2}+(4+6\nu)\,c_{1,4})\,\phi'\,\psi''
-((1+3\nu)\, c_{2,3}+(4+6\lambda)\,c_{4,1})\,\phi''\,\psi')
\end{array}
$$
\begin{enumerate}
\item If $(\lambda,\nu)=(0,0),$ then the constant $c_{1,4}$ is
arbitrary. On the other hand,
$$c_{0,5}=c_{5,0}, \quad c_{4,1}=-c_{1,4},\quad
c_{2,3}=6\,c_{4,1}, \quad c_{3,2}=-6\,c_{4,1}.$$ The constant
$c_{1,4}$ can be chosen in such a way that once adding the trivial
1-cocycle $L_XJ_5^{\lambda,\nu}$ above to our 1-cocycle
(\ref{co5}), the constant $\eta_2$ disappears completely.
Therefore, our 1-cocycle is generated by $\beta_1, \beta_2$ and
$\eta_3.$ Hence, the cohomology group is three-dimensional. The
result holds true for $(\lambda,\nu)=(0,-4), (-4,0).$

\item If $(\lambda,\nu)=(0,-2),$ then the constant
$c_{4,1}$ and $c_{0,5}$ are arbitrary. On the other hand,
$$c_{1, 4} = c_{4, 1}, \quad c_{2, 3} = 2\,c_{4, 1}, \quad
  c_{3, 2}= 2\, c_{4, 1}, \quad c_{5, 0}= \frac{1}{5}\,c_{4, 1}.$$
The constant $c_{1,4}$  and $c_{0,5}$ can be chosen in such a way
that once adding the trivial 1-cocycle $L_XJ_5^{\lambda,\nu}$
above to our 1-cocycle (\ref{co5}), the constant $\eta_1$ and
$\eta_4$ disappear completely. Therefore, our 1-cocycle is
generated by $\beta_3$ and $\eta_2.$ Hence, the cohomology group
is two-dimensional. The result holds true for
$(\lambda,\nu)=(-2,0)$ and $(-2,-2).$ \item If $(\lambda,
\nu)=(-\frac{3}{2}, -\frac{1}{2})$ then by Proposition
(\ref{gor}), the transvectant $J_5$ is not unique. A direct
computation proves that the constants $c_{5, 0}$ and $c_{2, 3}$
are arbitrary. On the other hand,
$$c_{0, 5} = \frac{c_{2, 3}}{10}, \quad
c_{1, 4}=\frac{c_{2, 3}}{2}, \quad c_{3, 2}= c_{2, 3}, \quad c_{4,
1}=5\,c_{5, 0}.
$$
The constant $c_{2,3}$ and $c_{5,0}$ can be chosen in such a way
such that once adding the trivial 1-cocycle $L_XJ_5^{\lambda,\nu}$
above to our 1-cocycle (\ref{co5}), the constant $\eta_1$ and
$\eta_4$ disappear completely. Therefore, our 1-cocycle is
generated by $\eta_2.$ Hence, the cohomology group is
one-dimensional. The result holds true for $(\lambda,
\nu)=(-\frac{1}{2}, -2), (-1,-1), (-1, -\frac{3}{2}),$ $(-1, -2),$
$(-\frac{1}{2}, -\frac{3}{2}),$ $(-\frac{3}{2}, -1),$
$(-\frac{3}{2}, -\frac{3}{2}), (-\frac{3}{2}, -2), (-2,
-\frac{1}{2}), (-2, -1), (-2, -\frac{3}{2}).$ \item If $\lambda$
and $\nu$ are not like above, then the transvectant $J_5$ is
unique by Proposition (\ref{gor}). Whatever the weights $\lambda$
and $\nu$ can take, the trivial 1-cocycle $ L_X J^{\lambda,\nu}_5$
is never identically zero. Therefore, one of the constants
$\gamma,$ $\alpha_1,$ $\alpha_2,$ $\beta_1,$ $\beta_2,$ $\beta_3,$
$\eta_1, \eta_2,\eta_3$ or $\eta_4$ can be eliminated by just
adding the trivial 1-cocycle $ L_X J^{\lambda,\nu}_5.$ Hence, the
cohomology group is two-dimensional.
\end{enumerate}
\subsection{The case when $k=6$}
If $\mu-\lambda-\nu=6,$ equation (\ref{sam2}) shows that the
1-cocycles that span the cohomology group of Theorem (\ref{main})
are of the form
\begin{eqnarray}
{\mathfrak
P}(X,\phi,\psi)&:=&\gamma\,X^{(7)}\,\phi\,\psi+X^{(6)}\,\left (
\alpha_{1}\,\phi'\,\psi +\alpha_{2}\,\phi\,\psi'\right
)+X^{(5)}\left
(\beta_1\,\phi''\,\psi+\beta_2\,\phi'\,\psi'+\beta_3\,\phi\,
\psi'' \right )\nonumber\\
&&+X^{(4)}\left
(\eta_1\,\phi^{(3)}\,\psi+\eta_2\,\phi''\,\psi'+\eta_3\,\phi'\,
\psi''+ \eta_4\,\phi\,\psi^{(3)}\right )\nonumber \\
&&+X'''\left
(\xi_1\,\phi^{(4)}\,\psi+\xi_2\,\phi^{(3)}\,\psi'+\xi_3\,\phi''\,
\psi''+ \xi_4\,\phi'\,\psi^{(3)} + \xi_5\,\phi\,\psi^{(4)} \right
) \label{co6}
\end{eqnarray}
The 1-cocycle condition is equivalent to the following system
$$
\begin{array}{ll}
14\gamma=-\lambda\,\beta_1-\nu\,\beta_3+\lambda
\,\xi_1+\nu\,\xi_5,
&14\gamma=-\lambda\,\alpha_1-\nu\,\alpha_2,\\
5\alpha_1=-(1+3\lambda)\,\eta_1+(1+4\lambda)\,\xi_1+\nu
(\xi_4-\eta_3),& 9\alpha_1=-(1+2\lambda)\,\beta_1-\nu \beta_2,
\\
5\alpha_2=-(1+3\nu)\,\eta_4+(1+4\nu)\,\xi_5+\lambda
(\xi_2-\eta_2), &
9\alpha_2=-(1+2\nu)\,\beta_3-\lambda\, \beta_2,\\
5\beta_1=-(3+3\lambda)\,\eta_1-\nu\,\eta_2,&
5\beta_3=-\lambda\,\eta_3-(3+3\nu)\,\eta_4\\
5\beta_2=-(1+2\lambda)\,\eta_2-(1+2\nu)\,\eta_3.&
2\eta_1=-(6+4\lambda)\, \xi_1 -\nu \, \xi_2,\\
2\eta_2=-(3+3\lambda)\, \xi_2 -(1+2\nu ) \xi_3,&
2\eta_3=-(1+2\lambda)\, \xi_3 -(3+3\nu)\xi_4,\\
2\eta_4=-(6+4\nu)\, \xi_5 -\lambda \, \xi_4,
\end{array}
$$
The space of solutions of the system above is four-dimensional for
$(\lambda, \nu)=(-\frac{5}{2}, -\frac{5}{2}),$ $(-\frac{5}{2},0),$
$(0,-\frac{5}{2}),$ $(0, \frac{-5\pm\sqrt{19}}{2}),$ $(
\frac{-5\pm\sqrt{19}}{2},0),$ $(\frac{-5-\sqrt{19}}{2},
\frac{-5+\sqrt{19}}{2}),$ $(\frac{-5+\sqrt{19}}{2},
\frac{-5-\sqrt{19}}{2}),$ and three-dimensional otherwise.
According to these values, let us study the triviality of the
1-cocycle (\ref{co6}). A direct computation proves that
\begin{gather*}
\nonumber L_XJ_5^{\lambda,\nu}(\phi,\psi)=(-\lambda
\,c_{6,0}-\nu\,
c_{0,6})X^{(7)}\,\phi\,\psi-((1+6\lambda)\,c_{6,0}+\nu\,c_{1,5})
X^{(6)}\phi'\,\psi\\\nonumber
\phantom{L_XJ_5}{-((1+6\nu)\,c_{0,6}+\lambda\,c_{5,1})
X^{(6)}\,\phi\,\psi'-
((6+15\lambda)\,c_{6,0}+\nu\, c_{2,4})\,X^{(5)}\,\phi''\,\psi}
\qquad\quad\\
\phantom{L_XJ_5}{-(((1+5\lambda)\,c_{5,1}+(1+5\nu)\,
c_{1,5})\,\phi'\,\psi'- ((6+15\nu)\,c_{0,6}+\lambda\,
c_{4,2})\,\phi\,\psi'')\,X^{(5)}}\quad\\
\phantom{L_XJ_5}{-((\nu\,c_{3,3}+(15+20\lambda)\,c_{6,0})\,\phi^{(3)}\,\psi
-((1+4\nu)\, c_{2,4}+(6+9\nu)\,c_{5,1})\,\phi''\,\psi'}\qquad\quad
\\
\phantom{L_XJ_5}{-((1+4\,\lambda)c_{4,2}+(6+9\nu)\,c_{1,5})\,\phi'\,\psi''
-(\lambda\,
c_{3,3}+(15+20\nu)\,c_{0,6})\,\phi\,\psi^{(3)}))\,X^{(4)}}
\\
\phantom{L_XJ_5}{-(((20+15\lambda)\,c_{6,0}+\nu\,c_{4,2})\,\phi^{(4)}\,\psi
-((10+10\lambda)\,
c_{5,1}+(1+3\nu)\,c_{3,3})\,\phi^{(3)}\,\psi'}\quad
\\
\phantom{L_XJ_5}{-((4+6\lambda)\,c_{4,2}+(4+6\nu)\,c_{2,4})\,\phi''\,\psi''
-((1+3\lambda)\, c_{3,3}+(10+10\nu)\,c_{1,5}) \,\phi'\,\psi^{(3)}}\\
\phantom{L_XJ_5}
{-(\lambda\,c_{2,4}+(20+15\nu)\,c_{0,6})\,\phi\,\psi^{(4)})\,X'''}\qquad\qquad
\qquad\qquad\qquad\qquad\qquad\qquad\quad
\end{gather*}
\begin{enumerate}
\item If $(\lambda,\nu)=(0,\frac{-5-\sqrt{19}}{2}),$ then the constant
$c_{5,1}$ is arbitrary. On the other hand,
$$c_{0, 6}= 0, \quad
  c_{1, 5}= -\frac{4}{45}(-35+8\sqrt{19})\,c_{5, 1},
  \quad c_{2, 4}= 2(-13+ 3\sqrt{19})\,c_{5, 1}, $$
  $$ c_{3, 3}= -\frac{4}{3}(-31+7\sqrt{19})\,c_{5, 1},
  \quad c_{4, 2}= \frac{10}{3}(-4+ \sqrt{19})\,c_{5, 1},
  \quad c_{6, 0}= \frac{1}{30}(5+ \sqrt{19})\, c_{5, 1}$$ The constant
$c_{5,1}$ can be chosen in such a way such that once adding the
trivial 1-cocycle $L_XJ_6^{\lambda,\nu}$ above to our 1-cocycle
(\ref{co6}), the constant $\xi_4$ disappears completely. Hence,
the cohomology group is three-dimensional. The result holds true
for $(\lambda,\nu)=(\frac{-5-\sqrt{19}}{2},0),
(\frac{-5+\sqrt{19}}{2},0), (0,\frac{-5+\sqrt{19}}{2}),
(\frac{-5-\sqrt{19}}{2},\frac{-5+\sqrt{19}}{2}),
(\frac{-5+\sqrt{19}}{2},\frac{-5-\sqrt{19}}{2}).$ \item If
$(\lambda,\nu)=(-\frac{5}{2},-\frac{5}{2}),$ then the constant
$c_{4,2}$ and $c_{1,5}$ are arbitrary. On the other hand,
$$c_{0, 6}= 0, \quad c_{1, 5}=-\frac{c_{4, 2}}{10}, \quad
    c_{2, 4}=c_{4, 2}, \quad c_{3, 3}= -2\, c_{4, 2}, \quad
    c_{5, 1}=-\frac{c_{4, 2}}{10}, \quad c_{6, 0}= 0.$$ The constant
$c_{4,2}$  and $c_{1,5}$ can be chosen in such a way that once
adding the trivial 1-cocycle $L_XJ_6^{\lambda,\nu}$ above to our
1-cocycle (\ref{co6}), two of the constants disappear completely.
Hence, the cohomology group is two-dimensional. The result holds
true for $(\lambda,\nu)=(-\frac{5}{2},0)$ , $(0,-\frac{5}{2}).$
\item If $(\lambda,\nu)=(-\frac{1}{2}, -2),$ then by Proposition
(\ref{gor}), the transvectant $J_6$ is not unique. A direct
computation proves that
$$
c_{1, 5}=6\,c_{0, 6}, \quad c_{2, 4}=\frac{5}{2}\,c_{5, 1}, \quad
c_{3, 3}=\frac{10}{3}\, c_{5, 1}, \quad c_{4, 2}= \frac{5}{2}
c_{5, 1}, \quad c_{6, 0}=\frac{1}{6}\,c_{5, 1}.
$$
The constant $c_{0,6}$ and $c_{5,1}$  can be chosen in such a way
that once adding the trivial 1-cocycle $L_XJ_6^{\lambda,\nu}$
above to our 1-cocycle (\ref{co6}), the constant $\xi_2$ and
$\xi_5$ disappear completely. Hence, the cohomology group is
one-dimensional. The result holds true for
$(\lambda,\nu)=(-2,-\frac{1}{2}),$ $(-\frac{1}{2},-\frac{5}{2}),$
$ (-\frac{5}{2}, -\frac{1}{2}),$ $(-1, -\frac{3}{2}),$
$(-\frac{3}{2}, -1),$ $(-1, -2),$ $(-2, -1),$ $(-\frac{5}{2},
-1),$ $ (-1,-\frac{5}{2}),$ $(-\frac{3}{2},-\frac{3}{2}),$ $(-2,
-\frac{3}{2}),$ $(-\frac{3}{2},-2),$ $(-\frac{3}{2},
-\frac{5}{2}),$ $(-\frac{5}{2}, -\frac{3}{2}), (-2,-2),$ $(-2,
-\frac{5}{2}),$ and $(-\frac{5}{2}, -2).$ \item If $\lambda$ and
$\nu$ are not like above, then the transvectant $J_6$ is unique by
Proposition (\ref{gor}). Whatever the weights $\lambda$ and $\nu$
can take, the trivial 1-cocycle $ L_X J^{\lambda,\nu}_6$ is never
identically zero. Therefore, one of the constants $\gamma,
\alpha_1,\ldots$  can be eliminated by just adding the trivial
1-cocycle $ L_X J^{\lambda,\nu}_6.$ Hence, the cohomology group is
two-dimensional.
\end{enumerate}
\subsection{The case when $k=7$}
The proof here is the same as in the previous section. We just
point out that the space of solutions of the 1-cocycle property
(\ref{sam2}) is four-dimensional for $(\lambda, \nu)=(0,-3), (-3,
0),$ and $(-3,-3);$ three-dimensional for $(\lambda,
\nu)=(-2,-2),$ $ (\lambda,0),$ $ (0, \nu),$ $
(\lambda,-6-\lambda),$ $ (-\frac{3}{2}, -3),$ $
(-3,-\frac{3}{2}),$ $ (-\frac{3}{2}, -\frac{3}{2}), $
$(-\frac{5}{2}, -1), (-1, -\frac{5}{2}),$
 $ (-\frac{5}{2},
-\frac{5}{2}),$ $(\sqrt{19}-1,\frac{-5-\sqrt{19}}{2}),$$
(\frac{-5-\sqrt{19}}{2},\sqrt{19}-1),$\\$
(-\sqrt{19}-1,\frac{-5+\sqrt{19}}{2}),$$
(\frac{-5+\sqrt{19}}{2},-\sqrt{19}-1),
(\frac{-5-\sqrt{19}}{2},\frac{-5-\sqrt{19}}{2}),
(\frac{-5+\sqrt{19}}{2},\frac{-5+\sqrt{19}}{2}),$ and
two-dimensional otherwise.
\subsection{The case when $k\geq 8$}
For $k\geq 8,$ the number of variables generating any 1-cocycle is
much smaller than the number of equations coming out from the
1-cocycle condition - for instance, for $k=8$ the $\sharp
\mbox{(Variables)}=28,$ while $\sharp \mbox{(Equations)}=33.$ For
generic $\lambda$ and $\nu,$ the number of equations will
generates a one-dimensional space, which give a unique cohomology
class. This cohomology class is indeed trivial because the
expression $L_XJ_k^{\lambda, \nu}$ is also a 1-cocycle.
\begin{rmk}{\rm
For $k\geq 8$ and for particular values of $\lambda$ and $\nu,$
the cohomology group $\mathrm{H}^1(\Vect(\mathbb{RP}^1),
\Sl(2,\mathbb{R}); \cD_{\lambda,\nu;\mu})$ may not be trivial. For
instance, for $k=8$ we have
$$
\mathrm{H}^1(\Vect(\mathbb{RP}^1), \Sl(2,\mathbb{R});
\cD_{0,-\frac{7}{2};\frac{9}{2}})\simeq \mathbb{R}.$$}
\end{rmk}
\subsection*{Acknowledgments}
I am grateful to V. Ovsienko for his constant supports.

\end{document}